\documentclass[12pt]{amsart}
\usepackage{amssymb, amsfonts, amsmath,pstricks,pst-plot}

\usepackage{fullpage}

\newcommand{\noi}{\noindent}

\parskip=\smallskipamount
\numberwithin{equation}{section}

\newtheorem{izr}{Theorem}[section]
\newtheorem{lm}[izr]{Lemma}
\newtheorem{pos}[izr]{Corollary}

\theoremstyle{definition}

\newtheorem{op}[izr]{Remark}

\newcommand{\Nn}{\mathbb{N}}
\newcommand{\Rr}{\mathbb{R}}
\newcommand{\Cc}{\mathbb{C}}

\newcommand{\cC}{\mathcal{C}}

\newcommand{\dok}{\noindent \it Proof. \rm }
\newcommand{\ra}{\rightarrow \:}

\newcommand{\ol}[1]{\overline{#1}}

\newcommand{\ve}{\varepsilon}
\newcommand{\psh}{plurisubharmonic }
\newcommand{\nbhd}{neighborhood }

\newcommand{\1}{\v{c}}

\newcommand{\supp}{{\rm supp }}
\newcommand{\inter}{{\mathrm {Int}}}

\newcommand{\Qed}{$\mbox{         }$ \hfill $\mbox{      } \spadesuit$ \\}

\def\inic{\setcounter{izr}{0}}

\begin{document}

\inic \pagestyle{plain} \setcounter{section}{0} \setcounter{tocdepth}{1}

\title{The generalized Oka-Grauert principle for 1-convex manifolds}
\author{Jasna Prezelj, Marko Slapar}
\address{Jasna Prezelj, Faculty of Mathematics and Physics, Department of Mathematics, University of
Ljubljana, Jadranska 21, SI-1000 Ljubljana, Slovenia; Marko Slapar, Faculty of Education, Department of Mathematics and Computing, University of
Ljubljana, Kardeljeva plo\v s\v cad 16, SI-1000 Ljubljana, Slovenia} \email{ jasna.prezelj@fmf.uni-lj.si, marko.slapar@pef.uni-lj.si}

\thanks{The authors work was supported  by the research program P1-0291 at the Slovenian Research Agency.}

\bigskip\bigskip\rm

%\begin{abstract}

%This paper presents a proof of the generalized Oka-Grauert principle for 1-convex manifolds: Every continuous mapping $X \ra Y$ from a 1-convex
%manifold $X$ to a complex manifold $Y$ which is already holomorphic on a neighborhood of the exceptional set is homotopic to a holomorphic one
%provided that either $Y$ satisfies CAP or we are free to change the complex structure on $X.$
%\end{abstract}

\maketitle
%\begin{center}
%\begin{minipage}{12cm}
%\tableofcontents
%\end{minipage}
%\end{center}
%\vspace{1cm}

\section{Introduction and the main theorem}
\inic

Let  $Y$ be a complex manifold satisfying convex approximation property (CAP) and let $X$ be an arbitrary Stein manifold. Then the Oka-Grauert (or homotopic) principle
holds for mappings $X \ra Y$ (\cite{f1}; the name Oka manifold has recently been suggested for such manifolds $Y$). This means that each homotopy class of mappings $X \ra Y$ admits a holomorphic representative. 

Manifolds satisfying CAP are in some sense `large'. As an example of a manifold failing to satisfy CAP consider the annulus $Y
= \{z \in \Cc, 1/2 < |z| < 2\}.$ Let $X = \{z, 1/3 < |z| < 3\}.$ There is plenty of continuous mappings from $X$ to $Y$ but no nontrivial holomorphic
ones. The reason is that $Y$ is `too small' for $X.$  If we are free to change the holomorphic structure on $X,$ then we can find for every continuous
mapping $f_0\! : X \ra Y$ another Stein structure $J_1$ on $X$ homotopic to the initial one and a holomorphic mapping $f_1\!: (X,J_1) \ra Y$ in
the same homotopy class as $f_0.$ In general, the change of structure depends on $Y$ as well as on $f_0.$ In the simple example above, the manifold $X$ is
homotopically equivalent to the unit circle $S^1 \subset X$ and we change the homotopic structure of $X$ just by squeezing it diffeotopically into a
small neighborhood of the unit circle. For a general Stein manifold $X$ we can proceed analogously, replacing $S^1$ by a suitable fattened $CW$-complex embedded in $X$ and
homotopically equivalent to $X$  to obtain the \\

\noi{\bf Generalized Oka-Grauert principle}: Every continuous mapping $X \ra Y$ from a Stein manifold $X$ to a complex manifold $Y$ is homotopic to a
holomorphic one provided that either $Y$ satisfies CAP or we are free to change the complex structure on $X$ (\cite{fs1}).\\ In addition
we can also require that the structure is fixed on a neighborhood of an analytic set $X_0 \subset X$ if the initial mapping is holomorphic on a
neighborhood of $X_0$ (\cite{fs2}).

It has recently been shown that if $X$ is 1-convex and $Y$ satisfies CAP, then the following version of the Oka-Grauert principle holds:\\

\noi{\bf Relative Oka-Grauert principle for mappings}: Every continuous mapping $X \ra Y$ from a 1-convex manifold $X$ to a complex manifold $Y$ that satisfies CAP which
is already  holomorphic on a neighborhood of the exceptional set is homotopic to a holomorphic map and the homotopy is fixed on the exceptional set
(compare \cite{lv}, \cite{p}).\\

Recall that a complex space $X$ is $1$-convex if it possesses a plurisubharmonic exhaustion function which is strictly \psh outside a compact set and
that the Remmert reductions of  $1$-convex spaces are  Stein. It is therefore possible to combine these two principles in
a\\

\noi{\bf Generalized Oka-Grauert principle for 1-convex manifolds}: Every continuous mapping $X \ra Y$ from a 1-convex manifold $X$ to a complex
manifold $Y,$ which is already  holomorphic on a neighborhood of the exceptional set is homotopic to a holomorphic map provided that either $Y$
satisfies CAP or we are free to change the complex structure on $X.$\\

\noi The following is the statement of the main theorem of this paper:

\begin{izr}[Generalized Oka-Grauert principle for 1-convex manifolds]\label{main thm}
  Let $(X, J_0)$ be a $1$-convex manifold of dimension at least $3$, $S$ its exceptional set,  $K \subset X$ a holomorphically convex compact
  subset of $X$ containing $S,$ $Y$ a complex manifold and $f_0\! : X \ra Y$ a
  continuous mapping which is holomorphic in a neighborhood of $K.$

  Then there exists a homotopy $f_t\! : X \ra Y$ and a homotopy $J_t$ of complex structures on $X$
  such that
  \begin{enumerate}
    \item[$(1)$] $f_t(x) = f_0(x)$ for $x \in S,$
    \item[$(2)$] $J_t = J_0$ on a neighborhood of $K$
    \item[$(3)$] $(X,J_1)$ is 1-convex with the exceptional set $S$,
    \item[$(4)$] mappings $f_t$ are $J_t$ holomorphic on a neighborhood  of $K$ and approximate $f_0$ on $K$
    as well as desired,
    \item[$(5)$]  $f_1$ is $J_1$ holomorphic on $X.$\\
  \end{enumerate}

\end{izr}
\vskip12pt

\section{Technicalities}
\inic

\subsection{Handle attaching in Stein category}

Let $(X,J)$ be a complex manifold with a complex structure $J$. A real immersed submanifold $i\!:\Sigma\rightarrow  X$ is {\it totally real} or {\it J-real} in $X$ if at every point $p\in \Sigma$ we have
$T_p i(\Sigma)\cap J(T_pi(\Sigma))={0}.$ The condition implies that dim$_\Rr$ $\Sigma \le$ dim$_\Cc$ $X$.

Let $W\Subset X$ be a relatively compact domain defined by $W=\{\rho<0\}$ where $\rho$ is a smooth real function defined in a neighborhood of $\partial
W$ and $d\rho\ne 0$ on $\partial W$. We say that $W$ (or $\partial W$) is {\it strongly pseudoconvex} or {\it J-convex}, if $dd^C\rho$ is a positive
form in a neighborhood of $\partial W$, meaning that $dd^C\rho(\nu,J\nu)>0$ for every $\nu\in TX|_{\partial W}$. So $dd^C\rho$ defines a metric on
$TX|_{\partial W}$ and in doing so, it also defines normal directions to $\partial W$. Now let  $W=\{\rho<0\}$ be a J-convex relatively compact domain in $(X,J)$,
$D=D_k\subset \Rr^k$ the closed unit ball with the boundary $S=S^{k-1}$. An {\it embedding (immersion) of a pair} $G\!:(D,S)\rightarrow (X\backslash
W,\partial W)$ is a smooth embedding (immersion) $G\!:D\rightarrow X\backslash W$, $G(S)=G(D)\cap \partial W$ with $G$ transverse to $\partial W$.
 We say that $G$ is {\it normal} to $\partial W$ if $dG$ maps normal vectors from $S=\partial D$ to normal vectors of $\partial W$.
$G(S)$ is {\it Legendrian} in $\partial W$ if $dG$ maps vectors tangent to $S$ into the {\it contact distribution} $T^{\Cc}W=T\partial W\cap JT\partial
W$ along $\partial W$.

The following lemma \ref{E1}, although not given in the above-stated manner, is proved in \cite{E}. A complete proof can also be found in \cite{fs1}. The main ingredients
are the Legendrization theorem of Gromov \cite{Gbook} and Duchamp \cite{Du} and the h-principle of Gromov for totally real submanifolds of complex
manifolds \cite{Gbook}. Lemma \ref{E1} functions just as well in an almost complex case.
\begin{lm}
\label{E1}
Let $W$ be an open, relatively compact J-convex domain with a smooth boundary in a complex manifold $(X,J)$ of the complex dimension
$n\ge 3$. For $0\le k\le n$ let $(D,S)$ be the closed unit disc in $\Rr^k$ with the boundary sphere $S$ and $G_0\!:(D,S)\to (X\backslash W,\partial W)$
a smooth embedding. Then there exists a regular homotopy of embeddings $G_t\!:(D,S)\to (X\backslash W,\partial W)$ $(0\le t\le 1)$ such that
\begin{enumerate}
  \item[$(1)$] $G_1$ is normal to $\partial W$,
  \item[$(2)$] $G_1(S)$ is Legendrian in $\partial W$,
  \item[$(3)$] $G_1(D)$ is totally real in $X$.
\end{enumerate}
If $\partial W$ is real analytic in a neighborhood of $G_0(S)$, then $G_1$ can also be made real analytic.
\end{lm}

\begin{op}
\label{E1dim2}
 In the complex dimension $2$, lemma \ref{E1} is also valid as stated if the attaching disc $D$ is one-dimensional. If the disc $D$ is two-dimensional,
 one cannot get the isotopy of embeddings but only a regular homotopy of immersions, so that the ending map is an embedding near the boundary of $D$
 but has special transverse double points in the interior of $D$. By special we mean that the double point is modeled by $\Rr^2\cup i\Rr^2\subset \Cc^2$.
\end{op}

Once we have a real analytic totally real disc $D$ attached normally (from the outside) along a Legendrian curve to a boundary $\partial W$ of a
strictly pseudoconvex domain $W\subset X$, we can use a holomorphic change of coordinates in a neighborhood of $\partial D$  coupled with a $\cC^1$
small real analytic deformation of the boundary $\partial W$ near $\partial D$ to get model situations of straight discs attached to a quadratic domain in $\Cc^n$. There we can use concrete functions to find strictly pseudoconvex neighborhoods of $W\cup D$ which preserve the topology of $W\cup D$. The
construction was first explained in \cite{E} and later also in \cite{FK}. More precisely, we have

\begin{lm} \label{E2} Let $W$ be an open, relatively compact J-convex domain in a complex manifold $(X,J)$. For $0\le k\le n$ let $(D,S)$
be the closed unit disc in $\Rr^k$ with the boundary sphere $S$ and $G_0\!:(D,S)\mapsto (X\backslash W,\partial W)$ a smooth totally real embedding with
$G(S)$ Legendrian. Let $d$ be any smooth Riemannian metric on $X$ and let $D_\ve$ be the $\ve$-neighborhood of $G(D)$ in $X$. Then for every $\ve> 0$
there exists a relatively compact strictly pseudoconvex Stein neighborhood $\tilde W$ of $W\cup G(D)$ with a smooth boundary completely contained in $W\cup D_\ve$ such that $W\cup G(D)$ is a
smooth deformation retract of $\tilde W$ and $W$ is Runge in $\tilde W$.
 \end{lm}

\subsection{Basics for $1$-convex manifolds}

The main problem with $1$-convex spaces is the lack of $1$-convex neighborhoods of graphs of holomorphic mappings $f\!: X \ra Y.$ However, if we remove
from the graph the zero set of a holomorphic function $g\!: X \ra \Cc$, which is zero on the exceptional set we get a Stein space and its graph consequently
has Stein neighborhoods. In addition, we do not want our Stein neighborhoods to be too `thin' in the $Y$-direction when approaching the graph of $f$
over $g^{-1}(0);$ we would like that the width in the $Y$-direction decreases at most polynomially. Such neighborhoods will be called {\it conic}
along the graph of $f$ over $g^{-1}(0).$ The existence of conic neighborhoods is given by

 \begin{izr}[Conic neighborhoods, \cite{p}, theorem 3.2] Let $X$ be a $1$-convex complex space with an exceptional set $S.$ Let
 $A \subset X$ be compact and holomorphically convex with $A \supset S$  and $Y$ a complex manifold. Let $f: X \ra Y$ be a continuous mapping,
  holomorphic on a $1$-convex neighborhood $U$ of $A$
and $g\! : X \ra \Cc$ a holomorphic function satisfying $g(S) = 0.$ Then for each  $1$-convex open set $U' \Subset {U}$ containing $A$ there exists a
Stein neighborhood of the graph $\Gamma f$ of $f$ over the set  $U'\setminus g^{-1}(0)$ in $X \times Y$  which is conic along the graph $\Gamma f $
over $g^{-1}(0).$
\end{izr}

In the theory of $1$-convex spaces there is a version of Cartan's theorem $B$ for relatively compact strictly pseudoconvex sets and there is also a version of
the Cartan's theorem A
\begin{izr}[Theorem A for $1$-convex spaces, \cite{p}] \label{thm A}
 Let $X$ be a $1$-convex space with an exceptional set $S,$  $U \Subset X$, $U$ open strictly pseudoconvex set
 containing $S,$
 ${\mathcal J} = {\mathcal J}(S)$ the ideal sheaf generated by
 the set $S$ and  ${\mathcal Q}$ a coherent sheaf on $X.$  There exists an $n_0 \in \Nn$ such that for  $n \geq n_0$
 the sheaf ${\mathcal Q} {\mathcal J}^n$ is locally generated by $\Gamma(U, {\mathcal Q} {\mathcal J}^n)$ on $U.$
\end{izr}

This theorem combined with the previously stated one  has many important consequences.

\begin{pos}[Existence of local sprays, \cite{p}, section 4.]\label{local sprays}
Let $X$ be a $1$-convex space, let
 $U' \Subset U \Subset X$ be  strictly pseudoconvex open sets containing the exceptional set $S$ of $X,$ and $f\! : X \ra Y$ a holomorphic mapping.
Then there exists a local spray on $U'$ fixing $S$ which dominates on $U' \setminus S,$ i.e. there exists a  holomorphic map $F\!: U' \times
B_n(0,\delta) \ra Y,$ such that $F(x,\cdot)\! : B_n(0,\delta) \ra Y,$ $D_tF(x,t)$ is surjective for $t = 0$ and $x \in U' \setminus S$,  $F(x,0) =
f(x)$ and   $F(x,t) = F(x,0)$ for all $x \in S.$
\end{pos}

\dok Denote by $VT(X \times Y)$ the kernel of the derivative  of the projection $X \times Y \ra Y$ and let $V$ be a conic Stein \nbhd of the graph of
$f$  over  $U \setminus g^{-1}(0)$ denoted by $\Gamma f|_{U \setminus g^{-1}(0)}$  for some holomorphic function $g\! : X \ra \Cc$ with zeroes on $S.$
By theorem $A$ for $1$-convex spaces (theorem \ref{thm A}) for each $k \in \Nn$ big enough there exist finitely many vector fields $h_1,\ldots,h_n$ of
the bundle $VT(X \times Y)|_{\Gamma f|_{U}}$ with zeroes of order (at least) $k$ on $\Gamma f|_{(g^{-1}(0))}$ generating $VT|_{\Gamma f_{(U \setminus
g^{-1}(0))}}.$ We extend these vector fields on $V$ and integrate them. Since $V$ is conic, the fields can be integrated for sufficiently small times
$t \leq t_0$ for all $x \in U' \setminus S$ (provided $k$ is big enough). Because of the zeroes on $\Gamma f|_{(g^{-1}(0))},$ we can extend the flows
of the fields on the graph of $f$ over $(g^{-1}(0))$ thus obtaining a map $F\! : U' \times B_n(0,\delta) \ra Z$ henceforth fulfilling all the requirements.\Qed

\subsection{Special pseudoconvex bumps}

In this subsection we construct special pseudoconvex bumps similar to those constructed in \cite{hl}. The main difference is that we have
a strictly pseudoconvex open set and a large disc attached to it.
\begin{lm}\label{bumps}
  Let $X$ be a $1$-convex manifold with the exceptional set $S$ and $\rho: X \ra \Rr$ a plurisubharmonic exhaustion function, which is strictly plurisubharmonic outside
  a holomorphically convex compact set $K \supset S.$ Let $0$ be a regular value for $\rho.$
  Let $A := \{\rho \leq 0\} \Subset X$ be a strictly pseudoconvex set containing $K$ and $D \subset X$ a Legendrian disc attached to the boundary of $A.$  We are going to construct a compact set $C \subset (A\setminus K) $ and a compact set $B \supset D$ such that
  \begin{itemize}
   \item[$(1)$]  $A \cup B$ is strictly pseudoconvex,
   \item[$(2)$] $A$ and $B$ have a bases of strictly pseudoconvex open neighborhoods $\{U_A\},$ $\{U_B\}$ respectively such that
         $\{U_A \cup U_B\}$ is a basis of strictly pseudoconvex open neighborhoods of $A \cup B;$ moreover, the sets $U_B$ are Stein,
   \item[$(3)$]  $A \cap B = C$ and $C$ is Runge in any of the neighborhoods $U_B,$
   \item[$(4)$]  separation property:  $(A \setminus B) \cup (B \setminus A) = \emptyset $ and $(\ol{U_A} \setminus \ol{U_B}) \cup (\ol{U_B} \setminus \ol{U_A})= \emptyset.$
   \end{itemize}
\end{lm}

\dok The assertions  derive from lemma \ref{E2}, \cite{hl}  and from the fact that we can reduce this situation to the model one with a straight disc attached to the quadratic domain in $\Cc^n.$ For the set $C$ we may take a (full and closed) $\ve$-torus $T_{\ve}$ around $\partial D$ intersected by $A.$ Using the  same methods as in \cite{hl} we can smooth the edges of $T_{\ve} \cap A$ to obtain a strictly pseudoconvex set $C$ with smooth boundary.  According to lemma \ref{E2} there is a $1$-convex neighborhood $U$ of $(\inter \,C) \cup D$ containing ${\inter} \, C.$ Let $B := \ol{U}.$ By definition
the set $C$ shares a piece of boundary with $A$ so the set $A \cup B$ is  strictly pseudoconvex  with the separation property.
Since $0$ is a regular value for $\rho$ we can produce a basis of neighborhoods $\{U_A\}$ and $\{U_B\}$ with the listed properties just by
taking a family of sublevelsets $\{\rho \leq \ve_n\}$ for some sequence $\ve_n \ra 0.$  \Qed
\vskip12pt

\section{Proof of the main theorem}
\inic

The proof follows the one presented in \cite{fs2} for the Stein case. The key idea in \cite{fs2} is the following. A suitably chosen strictly
plurisubharmonic exhaustion function for a Stein manifold $X$ gives an increasing sequence of strictly pseudoconvex sets $\{X_j\}$ such that there is
exactly one critical point in each of them. We move from $X_j$ to $X_{j+1}$ in two steps: firstly we have to cross the critical point by attaching a
suitable disc $D$ to $X_j$ and approximate our initial function with the one  holomorphic on a neighborhood of $D \cup X_j.$
Then, we have to `fatten' the union $D \cup X_j$ to get to $X_{j+1}.$  In the limit we get the desired structure on $X$ and a mapping, homotopic to the initial one.

Before we start to explain the differences in our approach let us mention that in practice we construct something that looks like a fattened CW-complex embedded in $X$ which is diffeotopically equivalent to $X$.
 %There is, of course, no difference in the `fattening' procedure.
 Following the above procedure in the case when $X$ is $1$-convex, we try to cross the critical point of the suitable exhaustion function by attaching a disc $D$ to
$X_{j-1}.$  Recall that the disc is `large', because we get it by pushing down the small disc given by Morse theory so that its boundary hits the boundary of $X_{j-1}$ in the proper way. If $X$ were Stein we could at this point use theorem 3.2 in \cite{f3} and get the desired  neighborhood of $X_{j-1} \cup D$ and  the approximation.
However,   one of the essential ingredients of the proof of theorem 3.2 in \cite{f3} is the fact that $X_j$ is
relatively compact and Stein. There seems to be (at least for authors) no obvious way of how to replace this with $1$-convexity; the Remmert reduction does not
help. Therefore we have to do the crossing via gluing, following the idea of Henkin and Leiterer. In their notation, we would like to attach to the set $A=X_j$ a `bone' $B \supset D$ such that $[A,B, A\cup B]$ is a (version of a) special pseudoconvex bump. In \cite{hl} the set $B$
 is a Stein neighborhood of the disc such that it shares a piece of boundary with $A$ and then $C$ is the intersection $A \cap B.$ In our case the situation is not local, therefore the only way to get Stein neighborhoods and approximations is by using theorem 3.2 in \cite{f3} and
 lemma \ref{E2}. Let $A,\ B,\ C$ be as in lemma \ref{bumps} and $f$ holomorphic on $\ol A$.
 %In the sequel we define the set $C$ to be a Stein subset of $A$ sharing a piece of boundary with $A.$ Let $B$ be a Stein set
 %given by lemma \ref{E2} such that $A \cup B$ is $1$-convex.
 By theorem 3.2 in \cite{f3} for the set $C$ we  obtain a thin `bone' $B'$ together with the desired holomorphic approximation $f'$ for $f;$ the better the approximation, the thinner the  bone $B'.$  We may assume that $A \cap B' = C.$ The gluing procedure for $f$ and $f'$ is performed by solving the $\ol{\partial}$-equation on the union $A \cup B$ so that the resulting
function is a small perturbation of the initial ones. The set where we solve $\ol{\partial}$-equation must be prescribed in advance since the norm of
the operator solving this equation depends on the set.

In the gluing procedure we solve the $\ol{\partial}$-equation on the union $A \cup B$ of the following type of forms. To explain the idea assume that
$f$ and $f'$ are functions. Consider the form
$$
   \omega = \ol{\partial} (f + \chi(f' - f))
$$
 using a suitable cut-off function $\chi$ with support of $d \chi$ on $C$ (see for example \cite{fp1} for details).
 Formally, this form is defined on $A \cup B'.$ But since its support
 is contained in $C,$ we can trivially extend it to $A \cup B.$
Now we solve $\ol{\partial}$-equation for the form $\omega$ on $A \cup B$ and get a function which is holomorphic over $A \cup B'.$ The sets $A \cup B$ and $A
\cup B'$ are homotopically equivalent and from here we can continue with the fattening procedure. The details are explained below.\\

Using the Remmert reduction choose a \psh exhaustion function $\rho\! : X \ra [0,\infty)$ such that $\rho^{-1}(0) = K$ and $\rho$ is strictly \psh
outside $K.$ Choose regular values $c_0<c_1<c_2\cdots$  for every $j$ such that there is exactly one critical point $p_j$ of $\rho$ contained in
$\rho^{-1}(c_{j-1},c_j)$ and let $k_j$ denote the Morse index of $p_j$ (recall that $k_j\le n$ for all $j$). The sublevelsets $X_j=\rho^{-1}[0,c_j)$
are topologically obtained by attaching a handle of index $k_j$ to the domain $X_{j-1}$. The handle in question is the thickening in $X$ of the stable
disc of $p_j$ of the gradient flow of $\rho$. Assume  that $f$ is already holomorphic in a neighborhood of $\overline X_{c_0}.$ We shall construct a
sequence of complex structures $J_j$ on $X$ and a sequence of maps $f_j\!:X\to Y$ such that

\begin{enumerate}
 \item for every $j$ the manifold $(X_j,J_j)$ is 1-convex with the exceptional set $S$,
 \item $J_j=J_{j-1}$ on a neighborhood of $\ol{X_{j-1}}$,
 \item $\ol {X_{j-1}}$ is Runge in $(X_j,J_j)$,
 \item $f_j$ is $J_j$ holomorphic on a neighborhood of $\overline {X_j},$
 \item $f_j$ and $f_{j-1}$ differ by less than $\ve 2^{-j-1}$ on $X_{j-1}$, and are equal on $S$,
 \item  $f_j$ is homotopic to $f_{j-1}$ by a homotopy that is $\ve 2^{-j-1}$ close to $f_{j-1}$ on $\ol X_{j-1}$.
\end{enumerate}

On the $j=0$ step we take $J_0$ to be the original complex structure on $X$ and $f_0$ to be the original function. By assumption, $f_0$ is holomorphic in a \nbhd
of $X_0$.

Let us now assume that the above conditions are met at the $j-1$ level. We set $A=X_{j-1}$ and $f_0=f_{j-1}$. Let $M$ be the stable manifold for the
critical point $p_j$ of the gradient flow of $\rho$, by using the metric associated with the positive Levi form $dd^c\rho$. The disc $D=M\cap (X\backslash
A)$ is a smooth disc attached transversely to the boundary of $A$. Using a small perturbation of our domain (or the function $\rho$), we can assume that
$\partial A$ is real analytic in a \nbhd of $\partial D$. We use lemma \ref{E1} with $W=A$ and $J=J_{j-1}$, to deform the disc $(D,\partial D)$ with a
small isotopy $(D_t,\partial D_t)\hookrightarrow (X\backslash A,\partial A)$ of pairs to a real analytic totally real disc, which we again call $D$,
$D\subset X\backslash A$, attached to $\partial A$ along a real analytic Legendrian curve.  Let $B$ and $C$ be as in lemma \ref{bumps}.  Recall that the boundary $\partial C$ agrees with $\partial A$ near $\partial D$ (see Figure 1).
\begin{figure}[ht]
\centering
\psset{unit=0.7cm, xunit=1.2, linestyle=solid, linewidth=0.5pt} %%
   \begin{pspicture}(-6,-5)(6,5)

    \pscustom[fillstyle=vlines, hatchcolor=gray]     % central part, colored
    {
    \pscurve[linewidth=0.5pt](0,-0.3)(2.5,-0.8)(3.2,-1.5)(4,0)(3.2,1.5)(2.5,0.8)(0,0.3)(-2.5,0.8)(-3.2,1.5)(-4,0)(-3.2,-1.5)(-2.5,-0.8)(0,-0.3)
    %(2.5,-0.8)(3,-1.5)(3.5,-1)(4,0)
    }

    \pscustom[fillstyle=hlines,hatchcolor=gray] {\pscurve[liftpen=1](5,4)(3,1.5)(2.5,0)(3,-1.5)(5,-4)                 % right hyperbola, yellow
                                                                                                                        %crosshatch
    }

    \pscustom[fillstyle=hlines,hatchcolor=gray] { \pscurve[liftpen=1](-5,4)(-3,1.5)(-2.5,0)(-3,-1.5)(-5,-4)             % left hyperbola, colored 5
    }

    \psline[linestyle=solid,linewidth=0.7pt](-2.5,0)(2.5,0)                                                          % the core disc

    \psdot[dotsize=5pt](-4.5,-2) \rput(-4.2,-1.9){$S$} %%%%%%%%Konec

% % %   NOTATION % %

  %%%%%%Zakomentirala jaz
   \rput(4.5,2){ $A$}
   %notation
   \rput(-4.5,2){ $A$} \psline[linewidth=0.2pt]{->}(-1,2)(-1,0.05) \rput(-1,2.3){ $D$}
   \psline[linewidth=0.2pt]{->}(1,-1.5)(1.3,-0.3) \psline[linewidth=0.2pt]{->}(1.3,-1.8)(3,-0.5) \psline[linewidth=0.2pt]{->}(0.7,-1.8)(-3,-0.5)
   \rput(1,-2){$B$} \rput(3.2,0){ $A\cap B$}

 \end{pspicture}
 \caption{A special pseudoconvex bump} \label{Fig1}
\end{figure}
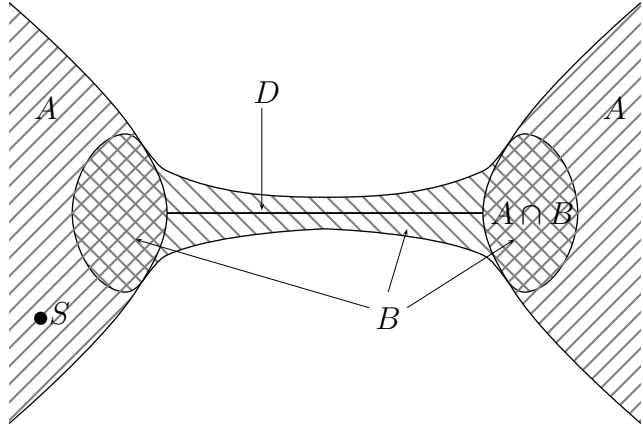
%%%%%%Konec

According to corollary \ref{local sprays} there exists a local spray $F$ on a neighborhood $U_A$ of $\ol A$ which dominates on $A\setminus S$ and keeps $S$ fixed.
According to theorem $3.2$ in (\cite{f3}) there exists an open `thin bone' $B'$ containing $\ol{C} \cup D,$   a Stein open set $U_{B'} \supset \ol{B'}$  and
 a map $G\! : U_{B'} \times B_n(0,\delta) \ra Y$ such that $G$ approximates $F$ on $C$ as well as desired.
There exists a `gluing' map $\gamma = (id_x,c)\!:
U_C \times B_n(0,\delta') \ra U_C \times B_n(0,\delta')$ such that $F = G \circ \gamma$ on $U_C \times B_n(0,\delta').$ The gluing map may be assumed to be close to the identity map. If it is not, then we can approximate it on $C$
 by $\tilde{\gamma} = (id_x, \tilde{c}) :U_B \times \Cc^n \ra U_B \times \Cc^n$ well enough so that $\tilde{c}$ is fiberwise invertible over $U_C$ and take $G \circ \tilde \gamma$ in the place of $G$ which gives $(id_x, \tilde{c}^{-1} \circ c)$ as a gluing map.

  Near the zero section the map $\gamma$ has a decomposition $\gamma = \beta
\circ \alpha^{-1}$ where $\alpha\! : U_A \times B_n(0,\delta') \ra U_A \times B_n(0,\delta'),$ and $\beta\! : U_B \times B_n(0,\delta') \ra U_B \times
B_n(0,\delta')$ are invertible on a \nbhd of the zero section over $U_C;$ $U_{B}$ is an open neighborhood of $\ol B$. Note that the decomposition is made
over the 1-convex open set $U_A \cup U_B.$ The desired decomposition of $\gamma$ is obtained from an implicit function theorem (Proposition 5.2 in
\cite{fp1}). In this theorem a solution of a $\ol{\partial}$-equation with uniform estimates on a \nbhd of strictly pseudoconvex set $A \cup B$ for
forms $\omega$ such that $\supp\ \omega \subset A \cap B$ is used where $A \cup B$ is a subset of a Stein manifold. Let us mention that this also works for Stein spaces.  In $1$-convex case the same can be
done with uniform estimates on the Remmert reduction of a neighborhood of $A \cup B$ since the supports of the forms we deal with do not intersect $S.$
In addition, we can require that the solutions of the $\ol{\partial}$ equation are zero on $S$ without spoiling the uniform estimates (see \cite{fp3}
for details). We have $F \circ \alpha (x,0) = G \circ \beta (x,0)$ for $x \in U_C.$ This defines a holomorphic map $f_1$ on a
  $A \cup B'$ homotopic to $f_0$ on $A \cup B'$. If $B'$ is thin enough, then according to lemma \ref{E2} there exists a strictly pseudoconvex open neighbourhood $\tilde W$
  of $A \cup D$ completely contained in $A\cup B'$. In addition, ${\tilde W} $ is diffeotopically equivalent to $X_j.$
  Outside ${\tilde W}$ the map $f_1$ can be glued to the map $f_0$ by a homotopy $f_t$ thus
  yielding a continuous map $f_1,$ holomorphic on a \nbhd of ${\tilde W},$  homotopic to $f_0$ and approximating $f_0$ on $A$. Since $F(x,t) = F(x,0)$
  for all $x \in S$ we also have $f_t(x) = f_0(x)$ on $S.$

%%%%Zakomentirala jaz
\begin{figure}[ht]
\centering
\psset{unit=0.7cm, xunit=1.2, linestyle=solid, linewidth=0.5pt} %%
\begin{pspicture}(-6,-4.5)(6,4.5)

%\pscustom[fillstyle=vlines,hatchcolor=gray]     % central part, colored
{
\pscurve[linewidth=0.8pt](0,-0.3)(2.4,-0.8)(4,-2.5)(5.5,0)(4,2.5)(2.4,0.8)(0,0.3)(-2.4,0.8)(-4,2.5)(-5.5,0)(-4,-2.5)(-2.4,-0.8)(0,-0.3)
}

%\pscustom[fillstyle=hlines,hatchcolor=gray]

{\pscurve[liftpen=1](5,4)(3,1.5)(2.5,0)(3,-1.5)(5,-4)                 % right hyperbola, yellow crosshatch
}

%\pscustom[fillstyle=hlines,hatchcolor=gray]
{ \pscurve[liftpen=1](-5,4)(-3,1.5)(-2.5,0)(-3,-1.5)(-5,-4)             % left hyperbola, colored
}

%\psline[linestyle=solid,linewidth=0.7pt](-2.5,0)(2.5,0)                                 % the core disc

\pscustom[fillstyle=hlines,hatchcolor=gray] { \pscurve(-3.4,2.1)(-3.7,0)(-3.4,-2.1)  % levi nosilec
\pscurve(-3.4,-2.1)(-3,-1.5)
\pscurve(-3,-1.5)(-3.2,0)(-3,1.5) \pscurve(-3,1.5)(-3.4,2.1) }

\pscustom[fillstyle=hlines,hatchcolor=gray] { \pscurve(3.4,2.1)(3.7,0)(3.4,-2.1)            % desni nosilec
\pscurve(3.4,-2.1)(3,-1.5)
\pscurve(3,-1.5)(3.2,0)(3,1.5) \pscurve(3,1.5)(3.4,2.1) }

\psdot[dotsize=5pt](-5.5,-2.5) \rput(-5.2,-2.4){$S$} %%%%%%%%Konec

% % %   NOTATION % %

\rput(5,3.5){ $A$}
 %notation
 \rput(-5,0){ $C$} \psline[linewidth=0.2pt]{->}(0.3,-1.8)(3.5,0)
\psline[linewidth=0.2pt]{->}(-0.3,-1.8)(-3.5,0) \rput(0,-2){ $\supp\ \omega$}
\end{pspicture}
\caption{Support of the form $\omega$} \label{Fig1a}
\end{figure}
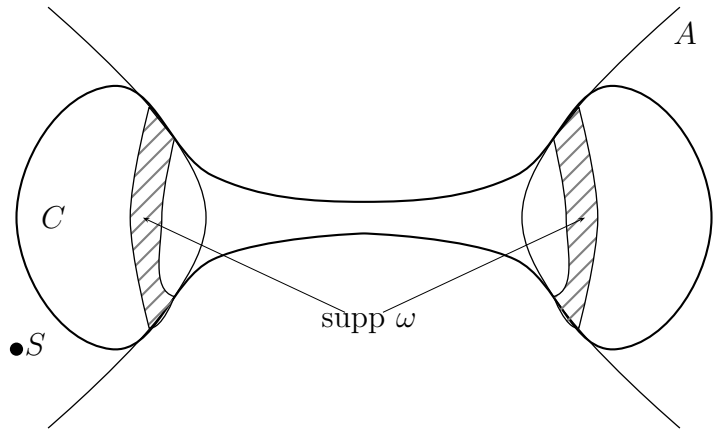
%%%%%%Konec
All that is needed now is to use a diffeotopy $H_t\!:X\to X$ with $H_0=Id$, $H_1(X_{j})=\tilde W$ and $H_t=Id$ on $X_{j-1}$. Such a diffeotopy exists
since $X_j$ is a noncritical extension of $\tilde W$. For every $t\in [0,1]$ let $J_t$ be the pull back of the almost complex structure $J_{j-1}$ by the
diffeomorphism $H_t$ and let $f_{1,t}$ be pull backs of $f_1$ above by $H_t$. We set $J_{j}=J_1$ and $f_{j}=f_{j,t}$.

 If there are only finitely many critical points for $\rho$, we have to  make another diffeotopy at the end in order to bring our whole manifold
 into the last sublevel set. If there are infinitely many critical points, the Runge condition ensures that the limiting complex structure on $X$
 does equip $X$  with a structure of 1-convex space with the singular set $S$. \Qed
\vskip12pt

\section{The case of dimension $2$}

In the complex dimension $2$ the situation is more complicated. The main difference is that lemma \ref{E1} is no longer valid in the complex dimension $2$ if
the attaching disc is also real $2$-dimensional. This was already noted by Eliashberg in \cite{E} and was subsequently justified using Seiberg-Witten
theory. For example, it is not possible to attach a totally real disc from the outside to the boundary of a unit ball in $\Cc^2$ along a Legendrian
curve. The obstruction essentially comes from the adjunction inequality for Stein surfaces, which for example, prohibits non null-homologous two spheres
from having self intersection number larger than $-2$. See for example \cite{Go1} or \cite{Nem}.

If a $2-$disc $D$ is attached to the boundary of a strictly pseudoconvex domain in a complex surface $X$ along a Legendrian curve $L$, then the Lai
indices \cite{lai} $I_{\pm}=e_{\pm}-h_{\pm}$ giving the difference of signed elliptic and hyperbolic complex points on $D$ are invariant under the
isotopy of embeddings that keeps the boundary fixed. The indices are calculated from the first Chern class of $X$ evaluated on $D$ and the relative self
intersection number of $D$. As long as $I_{\pm}$ both equal zero, one can find an isotopy of the disc $D$ through embeddings to a totally real disc by
keeping the boundary fixed. This is a result of Eliashberg and Harlamov \cite{EH}. The indices $I_{\pm}$ can always be increased by an isotopy of
embeddings that also moves the boundary $\partial D$ to a different Legendrian curve. It is, however, not possible to arbitrarily decrease the indices
$I_{\pm}$. Using remark \ref{E1dim2} following lemma \ref{E1} we can remedy this by introducing self-intersections (kinks) on the disc $D$ meaning
that by a regular homotopy of immersions we can make the disc $D$ to be an immersed totally real disc with its boundary circle an embedded Legendrian
curve. The immersed disc can be made to have only special double points so that the disc has tubular
Stein neighborhood basis.

To prove an analog of theorem \ref{main thm} in complex dimension $2$, we proceed as follows. As before, we decompose the manifold $(X,J)$ as an
increasing union of strictly pseudoconvex domains, so as to get to
 the larger domain: one either adds a handle of index $1$ or a handle of index $2$. If the critical point we have to pass to get to the next level
 has index $1$, there is no difference in the proof, since lemma \ref{E1} holds in this situation. Let us now  explain the difference when attaching
 an index $2$ critical point. As noted above, the disc $D$ attached to $A$ is in general not isotopic to a smooth disc $D'$ so that the union
 $C\cup D'$ has strictly pseudoconvex Stein neighborhoods, which was one of the essential ingredients in the proof in higher dimensions.
 The idea of how to fix this is Gompf's \cite{Go2} and is in a similar context explained in more detail in \cite{fs1}. First we find an isotopy of the disc $D$
 (a stable manifold of the critical point) so that it is Legendrian and normal at the boundary. If at the same time, we can also make the disc $D$
 totally real, we can just proceed as in the higher dimensional case.  If not, we add a sufficient number of positive standard kinks to $D$
 to get an immersed disc $D'$ so that $C\cup D'$ has thin tubular Stein strictly pseudoconvex neighborhoods.
 Adding a standard positive kink to $D$ means deleting a small disc in $D$ and along the boundary circle gluing back a small disc with exactly one positive transverse double point.

 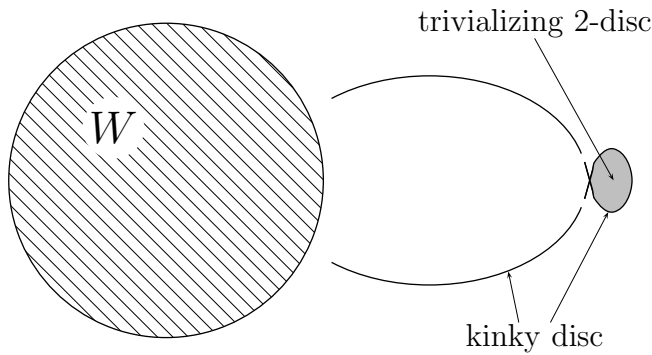
\begin{figure}[ht]
 \centering
\psset{unit=0.7cm, linewidth=0.5pt}

\begin{pspicture}(-5,-4)(5,4)
%%%%KOnec

\definecolor{myblue}{rgb}{0.66,0.78,1.00}

% % Leva elipsa %

\pscircle[fillstyle=vlines,hatchwidth=0.1pt,fillcolor=gray](-3,0){3} \pscircle[fillstyle=solid,fillcolor=white,linestyle=none](-4,1){0.6}
\rput(-4,1){\Large$W$} \psellipticarc(2,0)(3,2){10}{140} \psellipticarc(2,0)(3,2){220}{350} \psline(4.95,0.35)(5.05,0) \psline(4.95,-0.35)(5.05,0)

\psline(4.95,0.35)(5.14,-0.35) \psline(4.95,-0.35)(5.14,0.35)

\pscurve[fillstyle=solid,fillcolor=lightgray](5.05,0)(5.1,0.18)(5.14,0.35)(5.5,0.6)(5.85,0)(5.5,-0.6)(5.14,-0.35)(5.1,-0.18)(5.05,0)

\psline[linewidth=0.2pt]{->}(4.3,-2.7)(5.4,-0.6) \psline[linewidth=0.2pt]{->}(3.7,-2.7)(3.5,-1.73) \rput(4,-3){kinky disc} \rput(4,3){trivializing
$2$-disc} \psline[linewidth=0.2pt]{->}(4,2.7)(5.5,0)

\end{pspicture}
\caption{A kinky disc with a trivializing 2-disc} \label{Fig7}
\end{figure}

\noindent Although the relative homology class of $D'$ is the same as that of $D$, tubular neighborhoods of $D'$  are not diffeomorphic to the
thickening of $D$ because we have introduced extra generators in the $\pi_1$ group. To fix this, one adds a trivializing disc for each of the kinks,
and the thickening again becomes diffeomorphic to the thickening of $D$. Unfortunately, each of the trivializing discs also needs exactly one kink
to have thin tubular Stein neighborhoods. We therefore repeat the procedure. The limiting procedure (that is in this case necessarily infinite) gives an
(non smoothly) embedded disc $D''$ that agrees with $D$ near the boundary and has thin Stein neighborhoods homeomorphic but not necessarily
diffeomorphic to the tubular neighborhood of $D$. These limiting discs are called {\it Casson handles} and are an essential ingredient in the
classification theory of topological $4$-manifolds \cite{Fr,FQ}. Since for each critical point of index $2$ we may have to do infinitely many steps, we
use a variant of Cantor's diagonal process: first we make just one step on the new critical point and then go back to make one more step on each of the previous
critical points, before continuing to the next critical point. The analog of theorem \ref{main thm} is the following.

\begin{izr}[Generalized Oka-Grauert principle for 1-convex surfaces]\label{dim2}
  Let $(X, J)$ be a $1$-convex surface, $S$ its exceptional set,  $K \subset X$ a holomorphically convex compact subset of $X$ containing $S,$ $Y$ a
  complex manifold and $f\! : X \ra Y$ a
  continuous mapping which is holomorphic in a neighborhood of $K.$
  Then there exists a 1-convex surface $(X',J')$, a holomorphic map $f'\!:X'\to Y$ and an orientation preserving homeomorphism $h\!:X\to X$ so that
  \begin{enumerate}
    \item $h$ is holomorphic in a \nbhd of $K$ and $h(S)$ is the singular set for $X'$,
    \item $f'\circ h$ is homotopic to $f$ and $f'\circ h|_S=f|_S$,
    \item $f'\circ h=f$ approximates $f$ on $K$ as well as desired.
\end{enumerate}

\end{izr}

The point of the above theorem is that we may need to change the smooth structure on $X$.
The new results of Gompf (personal communication) indicate that the same theorem as for higher dimensional manifolds may hold in this case.

 \end{document}